\documentclass[review]{elsarticle}

\usepackage{lineno,hyperref}
\usepackage{amssymb, amsmath, amsthm}
\usepackage{multirow}
\usepackage{mathdots}
\newtheorem{theorem}{Theorem}[section]
\newtheorem{definition}{Definition}[section]
\newtheorem{lemma}{Lemma}[section]
\newtheorem{corollary}{Corollary}[section]

\newtheorem{Remark}{Remark}[section]
\newtheorem{Note}{Note}[section]
\newtheorem{Proposition}{Proposition}[section]
\modulolinenumbers[5]

\journal{Journal of \LaTeX\ Templates}
\newcommand{\least}{\let\cs=\@currsize\renewcommand{\baselinestretch}{.9}\tiny\CS}
\begin{document}
	\begin{frontmatter}
	\title{ Invariant bilinear forms under the operator  group of order $p^3$ with odd prime $p$ }
	
	\author[dam]{Dilchand Mahto \corref{cor1}}\ead{dilchandiitk@gmail.com}
	
	\author[dbm]{Jagmohan Tanti\corref{cor1}}\ead{jagmohan.t@gmail.com}
	
	\address[dam]{Department of Mathematics, Central University of Jharkhand, Ranchi, India}
	\address[dbm]{Department of Mathematics, Babasaheb Bhimrao Ambedkar University, Lucknow, India}
	
	
	
	
	
	
	
	\begin{abstract}
		In this paper we formulate the number of $n$ degree representations of a group of  order $p^3$ with $p$ an odd prime and the dimensions of  corresponding spaces of  invariant bilinear forms  over a field  $\mathbb{F}$ which contains a primitive $p^3$ root of unity. We explicitly discuss about the existence of a non-degenerate  invariant  bilinear form. 
	\end{abstract}
	
	\begin{keyword}
		\texttt  Bilinear  forms \sep Representation theory \sep  Vector spaces \sep  Direct sums \sep Semi direct product.
		\MSC[2010] 15A63\sep  11E04 \sep 06B15 \sep 15A03
	\end{keyword}
	
\end{frontmatter}


\section{Introduction}
Representation theory enables the study of a group  as operators on certain vector spaces and as an orthogonal group with respect to a corresponding bilinear form. Also since last several decades the search of non-degenerate invariant bilinear forms has remained of great interest. Such type of studies acquire an important place in quantum mechanics and other branches of physical sciences. \\
Let $G$ be a group and $\mathbb{V}$ a vector space over a field $\mathbb{F}$, then we have following.
\begin{definition}
	A homomorphism  $\rho$ : $G$ $\rightarrow$ GL($\mathbb{V}$) is called a representation of the group $G$.  $\mathbb{V}$ is also called a representing space of $G$. The dimension of $\mathbb{V}$ over $\mathbb{F}$ is called degree of the representation $\rho$.
\end{definition}
\begin{definition}
	A class function is a map  $f : G\rightarrow \mathbb{F}$ so that $f ( g ) =
	f ( h )$ if $ g $ is a conjugate of $h$ in $G$. 
	\label{def1.2}
\end{definition}
\begin{definition}
	A bilinear form on a finite dimensional vector space $\mathbb{V}(\mathbb{F})$ is said to be invariant under the representation $\rho$ of a finite group $G$ if 
	$$\mathbb{B}(\rho(g)x,\rho(g)y)=\mathbb{B}(x,y), \mbox{ $\forall$ g $\in$ $G$ and x,y $\in\mathbb{V}(\mathbb{F})$}.$$
\end{definition}
\noindent Let $\Xi$ denotes the space of bilinear forms on the vector space  $\mathbb{V}(\mathbb{F})$ over $\mathbb{F}$ and  $\mathcal{C}_{\mathbb{F}}$($G$) the set of all class functions on $G$.  
\begin{definition} The set of invariant bilinear forms under the representation $\rho$ is given by
	$$\Xi_{G} = \{\mathbb{B} \in \Xi \,\, |\, \mathbb{B}(\rho(g)x,\rho(g)y)=\mathbb{B}(x,y), \mbox{ $\forall$ $g$ $\in$ $G$ and x,y $\in$ $\mathbb{V}$}\}.$$
\end{definition}
\noindent Assume that $\mathbb{F}$ consists of a primitive $|G|^{th}$ root of unity. The  representation ($\rho$, $\mathbb{V}$) \cite{DCT} is irreducible of degree if and only if $\{0\}$ and  $\mathbb{V}$ are the only invariant  subspaces of $\mathbb{V}$ under the  representation $\rho$. The class function $\mathcal{C}_{\mathbb{F}}$($G$)   is a vector space over $\mathbb{F}$ with dimension $r$, where $r$ is the  number of  conjugacy classes of  $G$. By the Frobenius theorem (see pp 319, Theorem (5.9) \cite{Artin})  there are  $r$    irreducible representations $\rho_i$ (say), $1 \leq i \leq r$ of $G$ and $\chi_i$ (say) the corresponding characters of $\rho_i$.  
Also by Maschke's theorem ( see pp 316,  corollary (4.9) \cite{Artin}) every $n$ degree representation of $G$ can be written as a direct sum of copies of irreducuble representations. For  $\rho$ = $\oplus_{i=1}^{r}k_{i}\rho_{i}$ an $n$ degree representation of $G$, the coefficient of $\rho_{i}$ is $k_{i}$, $1 \leq i \leq r$, so that $\sum_{i=1}^{r} d_{i}k_{i}=n, $ and $\sum_{i=1}^{r} d_{i}^{2}= |G|, $ where $d_{i}$ is the degree of $\rho_{i}$  and $d_{i}||G|$ with  $d_{j}\geq d_{i}$ when $j>i$. It is already well understood in the literature that  the   invariant space $\Xi_{G}$ under $\rho$ can be expressed by the set  $\Xi_{G}'$ = $\{X \in \mathbb{M}_{n}(\mathbb{F})\,|\, C_{\rho(g)}^{t}X C_{\rho(g)} = X, \forall g \in G  \}$ with respect to an ordered basis $\underline{e}$ of $\mathbb{V}(\mathbb{F})$, where  $\mathbb{M}_{n}(\mathbb{F})$ is the set of square matrices of order $n$ with entries from $\mathbb{F}$ and  $C_{\rho(g)}$ is the matrix representation of the linear transformation $\rho(g)$ with respect to $\underline{e}$.\\

Here we consider $G$ to be a group of order $p^3$ with $p$ an odd prime, $\mathbb{F}$ a field with char($\mathbb{F}$) $\neq p$, which consists of a primitive $p^3$th root of unity  and   ($\rho$, $\mathbb{V}$)  an $n$ degree representation of $G$ over  $\mathbb{F}$. Then the corresponding  set $\Xi_{G}$ of  invariant bilinear forms on $\mathbb{V}$ under $\rho$, forms a subspace of  $\Xi$. In this paper our investigation is about the following questions.

\vskip1mm

\noindent\textbf{Question.} How many $n$ degree  representations (upto isomosphism) of  $G$  can be  there ? What is the dimension of  $\Xi_{G}$  for every $n$ degree representation ? What are the necessary and sufficient conditions for the existence of a non-degenerate invariant bilinear form.

The primary focus is on the existence of a non-degenerate invariant bilinear forms. Over complex numbers it has been seen with positive findings, as an evidence we present here one.

\vskip1mm

It is well known that every maximal (proper) subgroup of $G$ has index $p$ and is normal (As finite $p$ groups are nilpotent and any proper subgroup of a nilpotent group is properly contained in its normaliser). Thus there are epimorphisms from $G$ to the cyclic group $C$ of order $p$. 

Fix a generator $c$ of $C$ and $1\not=\zeta$ a primitive $p$th -root of unity. Let $U$ and $V$ be the one-dimensional repre-
sentations of $C$ on which $c$ acts respectively by $\zeta$ and $\bar\zeta$. We claim that $U\oplus V$ admits a $C$-invariant
non-degenerate bilinear form. Via some epimorphism to $C$ one can pullback these representations and the forms to $G$.

To prove the claim let us fix the vectors $0\not=u\in U$ and $0\not=v\in V$. Using these we define a bilinear form $B$ on $U\oplus V$ as follows: $B(u,u)=0=B(v,v)$, $B(u,v)=1=B(v,u)$ so that $B(\lambda u+\mu v, \lambda' u+\mu'v)=\lambda\mu'+\lambda'\mu$.

Now we may easily check the $C$ invariance as follows: $B(c(\lambda u+\mu v), c(\lambda' u+\mu' v))=B(\zeta(\lambda u+\mu v), \bar\zeta(\lambda' u+\mu' v))=\zeta\bar\zeta\lambda\mu'+\lambda'\mu=B(\lambda u+\mu v, \lambda' u+\mu'v)$.

The questions in concern have been studied in the literature in several distinct contexts.  Gongopadhyay and Kulkarni \cite{Gong2} studied the  existence of $T$-invariant non-degenerate symmetric (resp. skew-symmetric) bilinear forms. Kulkarni and Tanti \cite{Kulk} formulated the  dimension of space of  T-invariant bilinear forms. Mahto and Tanti \cite{DCT} formulated  the dimensions of invariant spaces and explicitly discussed about the existence of  the non-degenerate invariant bilinear forms under $n$ degree representations of a group of order 8. Sergeichuk \cite{VVSE} studied systems of forms
and linear mappings by associating with them self-adjoint representations of a category with involution. Frobenius \cite{GFRO} proved that every endomorphism of a
finite dimensional vector space V is self-adjoint for at least one non-degenerate symmetric
bilinear form on V. Later, Stenzel \cite{HSTE} determined when an endomorphism could be skew-
selfadjoint for a non-degenerate quadratic form, or self-adjoint or skew-self adjoint for a
symplectic form on complex vector spaces. However his results were later generalized to
an arbitrary field \cite{RGTL}. Pazzis \cite{CSPA} tackled the case of the automorphisms of a finite dimensional vector space that are orthogonal (resp. symplectic) for at least one non-degenerate
quadratic form (resp. symplectic form) over an arbitrary field of characteristics 2.

In this paper we investigate about the  counting of $n$ degree representations of a group of order $p^3$ with $p\geq3$ a prime, over a field $\mathbb{F}$ which consists of a primitive $p^3$th root of unity,  dimensions of their corresponding  spaces of invariant bilinear forms and establish a characterization criteria for existence of a non-degenerate invariant  bilinear form. Our  investigations  are stated in the following three main theorems. \\
\begin{theorem}
	The number of  n degree representations (upto isomorphism) of a group $G$ of order $p^{3}$, with $p$ an odd prime  is $\binom{n+p^{3}-1}{p^{3}-1}$ when $G$ is abelian and  $\sum_{\mu=0}^{\lfloor \frac{n}{p}\rfloor} \binom{\mu+ p-2}{p-2}\binom{n-\mu p+p^{2}-1}{p^{2}-1}$ otherwise.
	\label{theorem1.1}
\end{theorem}

\begin{theorem}
	The space $\Xi_{G}$ of invariant bilinear forms  of a group $G$  of order $p^3$ ($p$ an odd prime), under an $n$ degree representation 
	$(\rho ,\mathbb{V}(\mathbb{F}))$ is isomorphic to  the direct sum of  the subspaces 
	$\mathbb{W}_{(i,j)}$, $(i,j) \in A_G$  of $\mathbb{M}_{n}(\mathbb{F})$, i.e.,
	$\Xi_{G}'$ = $\bigoplus_{(i,j) \in A_G} \mathbb{W}_{(i,j)}$,
	where $A_G=\{(i,j)\,|\,\, \rho_i$ and $\rho_j$ are dual to each other\}  and  for every $(i,j)\in A_G$,  $\mathbb{W}_{(i,j)}$ = \{ $ X \in \mathbb{M}_{n}(\mathbb{F})\, | \, X_{d_{i}k_{i} \times d_jk_j}^{ij}=  C_{k_{i}\rho_{i}(g)}^{t}X_{d_{i}k_{i} \times d_jk_j}^{ij}C_{k_{j}\rho_{j}(g)}$, $\forall g \in G$ and rest blocks are zeros \}. Also for $(i,j)\in A_G$, the  dimension of $\mathbb{W}_{(i,j)}= k_ik_j.$
	
	
	.
	
	\label{theorem1.2}
\end{theorem}

\begin{theorem}
	
	If $G$ is group of order $p^3$, with $p$ an odd prime,  then an $n$ degree representation of $G$ 
	consists  of a non-degenerate invariant bilinear form if and only if every  irreducible representation and its dual  have the same multiplicity.\\
	
	\label{theorem1.3}
\end{theorem}

Thus we are able to give the answers to the questions in concern. Here for $n\in\mathbb{N}$ Theorem \ref{theorem1.1} counts all $n$ degree representations, Theorem \ref{theorem1.2} computes the dimension of the space of invariant bilinear forms and Theorem \ref{theorem1.3} characterises those $n$ degree representations, each of which admits a non-degenerate invariant bilinear form for a group of order $p^3$ with $p$ an odd prime over the field $\mathbb{F}$.

\begin{Remark}
	Thus we get the necessary and sufficient condition for the existence of a non-degenerate invariant bilinear form under an n degree representation.
\end{Remark}
\section{Preliminaries }
The classification of groups of order  $p^3$, with $p$ an odd prime has been well understood in the literature. Due to the structure theorem of finite abelian groups, there are only three abelian groups (upto isomorphism) of this order viz, $\mathbb{Z}_{p3} $, $\mathbb{Z}_{p^2} \times \mathbb{Z}_p$  and $\mathbb{Z}_p \times \mathbb{Z}_p \times \mathbb{Z}_p$. Amongst   non-abelian groups of this order Hisenberg group \cite{KCON} is well known and  named after a German theoretical physicist Werner Heisenberg. In this group  every non identity element  is of   order $p$. The elements of this group are usualy seen in the form of $3 \times 3$ upper triangular  matrices whose diagonal entries consist of 1 and other three entries are chosen from the finite field $\mathbb{Z}_p$. If at all  there exists any other non-abelian group of this order then it must have  a non identity element of  order  $p^2$. Let us consider the $2 \times 2$ upper triangular  matrices with $a_{11}=1+pm, (m \in \mathbb{Z}_p)$, $a_{12}=a \in \mathbb{Z}_{p^2}$ and $a_{22}=1$. Here the element with entries $a_{11}=a_{12}=a_{22}=1$ has order $p^2$ making it  non-isomorphic to the  Heisenberg group. We denote this group by $G_p$. Thus upto isomorphism there are  five  groups of order $p^3$ with an odd prime $p$ \cite{KCON}. For an abelian group of order $p^3$, there are $p^3$ number of  irreducible representations each having degree $1$   and for non-abelian cases, the number of  trivial conjugacy calsses is $|Z(G)| = p$. To find a non-trivial conjugacy class we refer to the theory of group action   and class equation  $$\displaystyle|G|= |Z(G)|+\sum_{\substack{g \notin Z(G)\, \mbox{varying over distinct conjugacy classes}}} |C_g|$$ with $|C_g| =  \frac{|G|}{|C(g)|}$, where $C_g$ and $C(g)$ are  the conjugacy class and the centralizer respectively of $g$ in $G$. If  $g \notin Z(G)$  then $|C(g)|=p^2$. Therefore there  are $p^2-1$ non trivial conjugacy classes of order $p$. Thus total number of conjugacy classes for a non-abelian group is $r =p^2 -1+p$, which is same as the  number of irreducible representations   with degree $d_i$ since  $d_i | |G|$ and $\sum_{i=1}^{r} d_i^2 =|G|$, therefore $d_i =1$ or $p$. Thus there are $p^2$ representations of degree $1$ and $p-1$ representations of degree $p$ for a non-abelian  group. We here formulate every  irreducible representation $\rho_i$ in such a way that the  entries of  $C_{\rho_i(g)}$ are either $0$ or $p^3$th primitive roots of unity. \\

\begin{definition}
	The  character  of $\rho$ is a  function    $\chi$ : $G$ $\rightarrow$ $\mathbb{F}$,  $\chi(g)=$ tr($\rho(g)$)  and is also called character of the group $G$.
\end{definition}
\begin{theorem}
	(Maschke's Theorem): If char($\mathbb{F}$) does not divide $|G |$, then every representation of  $G$ is a direct sum of irreducible representations.
\end{theorem} 

\begin{proof} 
	See p -316,   corollary (4.9)  \cite{Artin}.
\end{proof}

\begin{theorem}
	Two representations ($\rho, \mathbb{V(\mathbb{F})}$) and ($\rho', \mathbb{V(\mathbb{F})}$) of $G$ are  isomorphic iff their character values are same i.e, $\chi(g)=\chi'(g)$ for all g $\in G$.
	\label{theorem2.2}
\end{theorem}

\begin{proof}
	See  p -319,  corollary ( 5.13) \cite{Artin}.
\end{proof}
\section{Irreducible representation (irrep.) of group  of order $p^3$ with an odd prime $p$. }
\label{section3}
\noindent In this section $G$ is a group of order $p^3$ with $p$ an odd prime, $(\rho_i, \mathbb{W}_{d_i}(\mathbb{F})) $  stands for an  irreducible representation $\rho_i$ of degree $d_i$ of $G$  over a  field $\mathbb{F}$ with char($\mathbb{F}$) $\neq p$, which consists of $\omega \in \mathbb{F}$, a  primitive $p^3$th  root of unity.\\
Let $\sigma_s= \rho_{p^2+s}$, $1 \leq s \leq p-1$ denote the   irreducible representations of  degree $p$ when $G$ is non-abelian.  Since  $\sigma_s$ is a homomorphism from $G$ to $GL(\mathbb{W}_p) \cong GL(p, \mathbb{F})$, by the fundamental  theorem of homomorphism $\frac{G}{Ker(\sigma_s)} \cong $ $\sigma_s(G)$  and the  possible value of $|Ker(\sigma_s)|$ is $ 1 $ or $ p$. If $|Ker(\sigma_s)|= p$ then $(\chi_s ,\chi_s) \geq 1$, therefore $(\chi_s ,\chi_s) = 1$ only when $g \notin Ker(\sigma_s)$, so we have  $\chi_s(g)=0$. Also we have trivial character $\chi_i(g)=1, \forall g \in G$ and $(\chi_i ,\chi_s) = \frac{p^2}{p^3}\neq 0$, which  fails  the orthonormality property of   the irreducible characters. Thus $|Ker(\sigma_s)|= 1 $ and hence $\sigma_s(G)$ $\cong \frac{G}{\{e_G\}}  $ which  is isomorphic to a non-abelian group of order $p^3$. Now $\sigma_s(G)$ has  subgroups $H$ and $K$ of orders $p$ and $p^2$ respectively [see p-132, \cite{Dummit}, exercise-29]. Since $K$ is a maximal subgroup of $G$ so $K$  must be a normal subgroup and there exists a subgroup $H_p$ of order $p$ which is not   normal (not contained in $Z(G)$)[see p-188, \cite{Dummit}, Theorem 1], thus we have, $ \sigma_s(G) = KH_p$ and every element of $ \sigma_s(G) $ can be expressed uniquely in the form of $kh$ for some  $k \in K$ and $h \in H_p$ (this uniqueness follows from the condition $K \cap H_p = \{\sigma_s(e_G)\}$), therefore   $ \sigma_s(Heis(\mathbb{Z}_p)) \cong (\mathbb{Z}_{p} \times \mathbb{Z}_{p}) \rtimes \mathbb{Z}_{p}$ and  $ \sigma_s(G_p) \cong \mathbb{Z}_{p^2}  \rtimes \mathbb{Z}_{p}$ . As $|Z(G)|=p$, we can choose a  subgroup of order $p$ from $GL(p, \mathbb{F})$ and say that it is the Image of $Z(G)$ (subset of a normal subgroup of order  $p^2$  of $G$) denoted by $Im(Z(G))$ = $ \{ \omega^{sp^2}I_p  | \mbox{ $ 1 \leq s \leq p$ } \}$  under the irreducible representation $\sigma_s$ (since center elements commute and are  scalar matrices). Thus $Im(Z(G)) (\subsetneqq \sigma_s(K))\cong \mathbb{Z}_{p}$, each non-identity  element of  $Z(G)$ have $p-1$ choices in $\mathbb{Z}_{p}$ under $\sigma_s$ and rest  $p^3 -p$  elements  of $G$ map to rest $p^3 -p$  elements  of $\sigma_s(G)$. We decide all $p$ degree representations  by  the elements of   $Z(G)$  and   elements of $G -  Z(G)$ by  mapping to the set $\sigma_s(G)-ImZ(G)$  $\subseteq$ $GL(p, \mathbb{F})$, which consists of those elements whose trace is zero and order of every element is either $p$ or $p^2$. We will depict all irreducible   representations  of $G$ in the  next subsections. \\

\subsection{Heisenberg group.}
\begin{center}
	$G=$	$Heis(\mathbb{Z}_{p})$ = $\left\{ ((\alpha,\beta),\gamma)= \begin{bmatrix}
		1&\alpha&\beta \\
		0 &1&\gamma \\
		0&0&1
	\end{bmatrix}\,|\,\alpha,\beta,\gamma \in \mathbb{Z}_{p} \right\}$.
	\label{subs3.1}
\end{center}
Another presentation of Heisenberg group is (\cite{Dummit}, pp 179) 
\begin{center}
	$G= \langle x,a \, | \, x^p=a^p=1,\, x(a^{-1}xa)x^{-1}=a^{-1}xa,\, a(xax^{-1})a^{-1}=xax^{-1}\rangle$.
\end{center}

For the Heisenberg group, center 	$Z(G) =\langle xax^{-1}a^{-1}\rangle=\langle((0,1),0)\rangle$. Each of the irreducible representations $\sigma_s$, $1\leq s\leq p-1$ of degree $p$ maps the center of $G$ to the center of $GL(p, \mathbb{F})$, i. e.,  if $z \in Z(G)$ then $\sigma_s(z)$ is a scalar matrix say $c_sI_p$, $c_s \in \mathbb{F}$ and since order of $z$ is $p$ so $c_s^p=1$. Hence $c_s \in \{1,\omega^{p^2}= \epsilon, \epsilon^{2}, \cdots \cdots ,  \epsilon^{(p-1)}  \}$. Thus each of the $p-1$ irreducible representations of degree $p$ maps  $Z(G)$ into the $Z(GL(p, \mathbb{F}))$ even maps to the   $Z(\sigma_s(G))$, it has been recorded  in the following table. \\

\begin{table}[h]
	\centering
	\caption{All $p$ degree irreducible representations of $Heis(\mathbb{Z}_{p})$}
	\begin{tabular}{|c|c|c|c|c|c|c|c|}
		\hline
		All irrep. $\rightarrow$ &$\sigma_{2\eta-1}$&$\sigma_{2\eta}$\\
		\hline
		Runing variable $\rightarrow$	& \multicolumn{2}{c|}{1 $\leq \eta \leq \frac{p-1}{2}, $} \\
	
		\hline
		$g \in Z(G)$ & $\epsilon^{\eta } I_p$&$\epsilon^{(p-\eta)} I_p$\\
		\hline
		$g \notin Z(G)$& See Note \ref{note3.1}&See Note \ref{note3.1}\\
		\hline
		
		Dual irrep.	& \multicolumn{2}{c|}{$\sigma_{2\eta-1}^*=\sigma_{2\eta}$} \\
		\hline
		$\#$ irrep.	& \multicolumn{2}{c|}{$p-1$} \\
		\hline
		
	\end{tabular}.
	\label{table1}
\end{table}
\begin{Note}
	For  the $p$ degree representations  $\sigma_s$ of Heisenberg group $G$, if $xax^{-1}a^{-1} \in Z(G)$ then $\sigma_s(xax^{-1}a^{-1})= \rho_{p^2+s}(xax^{-1}a^{-1})=\epsilon^{m}I_p$, for some $m$, $1\leq  m < p$ and 
	the elements  of $G - Z(G)$ get mapped bijectively to the following set
	\begin{center}
		$\sigma_s(G)-\sigma_s\Big(Z(G)\Big)$  = $\{A_{u} \in GL(p, \mathbb{F})\, |\, \mbox{$Tr(A_{u}) = 0$ and $A_{u}^{p}= I_{p} $ for $1 \leq u \leq p^3 -p$}\}$.
	\end{center}
	\label{note3.1}
\end{Note}
Thus the $p$ degree irreducible  representations $\sigma_s$ can be expressed as below 
$$ \sigma_{2 \eta-1}(x)	=\begin{bmatrix}
	0&1&0&0&\cdots&0\\
	0&0&1&0&\cdots&0\\
	0&0&0&1&\cdots&0\\
	0&0&0&0&\cdots&0\\
	\vdots& \vdots&\vdots&\vdots&\ddots&\vdots\\
	1&0&0&0&\cdots&0
\end{bmatrix}^{\eta}\,  and \,\sigma_{2 \eta-1}(a)	=\begin{bmatrix}
	\epsilon^{\eta}&0&0&0&\cdots&0\\
	0&\epsilon^{\eta +1}&0&0&\cdots&0\\
	0&0&\epsilon^{\eta +2}&0&\cdots&0\\
	0&0&0&\epsilon^{\eta +3}&\cdots&0\\
	\vdots& \vdots&\vdots&\vdots&\ddots&\vdots\\
	0&0&0&0&\cdots&\epsilon^{\eta +p-1}
\end{bmatrix},$$

$$ \sigma_{2 \eta}(x)	=\begin{bmatrix}
	0&1&0&0&\cdots&0\\
	0&0&1&0&\cdots&0\\
	0&0&0&1&\cdots&0\\
	0&0&0&0&\cdots&0\\
	\vdots& \vdots&\vdots&\vdots&\ddots&\vdots\\
	1&0&0&0&\cdots&0
\end{bmatrix}^{-\eta} \,\, and \,\,\sigma_{2 \eta}(a)	=\begin{bmatrix}
	\epsilon^{p+1-\eta}&0&0&0&\cdots&0\\
	0&\epsilon^{p+2-\eta }&0&0&\cdots&0\\
	0&0&\epsilon^{p+3-\eta }&0&\cdots&0\\
	0&0&0&\epsilon^{p+4-\eta}&\cdots&0\\
	\vdots& \vdots&\vdots&\vdots&\ddots&\vdots\\
	0&0&0&0&\cdots&\epsilon^{p-\eta }
\end{bmatrix}.$$ 	

Non-abelian  groups of order $p^3$ 
have $p^2$ representations of degree  1. Let $\sigma_{(s-1,t-1)}$, $1  \leq  s,t \leq p$, denote representations of degree 1.  Here we present all 1 degree representations  for the  Heisenberg group. \\


\begin{table}[h]
	\centering
	\caption{All irreducible representations of degree 1 for the  Heisenberg group.}
	
	\begin{tabular}{|c|c|c|c|c|c|}
		
		\hline
		
		All irrep. $\rightarrow$&$\rho_{(t-1)p+2s-2}=$&$\rho_{(t-1)p+2s-1}=$&$\rho_{(t-1)p+2s-1}=$&$\rho_{(t-1)p+2s}=$\\
		& $\sigma_{(s-1,t-1)}$& $\sigma_{(p-s+1,p-t+1)}$ & $\sigma_{(s-1,t-1)}$& $\sigma_{(p-s+1,p-t+1)}$\\
		\hline
		
		Runing Variable	$\rightarrow$ &\multicolumn{2}{c|}{$1 \leq s \leq \frac{p+1}{2}, t=1$}& \multicolumn{2}{c|}{$1 \leq s \leq p,\,\, 2 \leq t \leq \frac{p+1}{2}$} \\
		\hline
		$\begin{bmatrix}
			1&0&\beta \\
			0 &1&0 \\
			0&0&1
		\end{bmatrix}$& $\epsilon^{p}=1 $& $1 $ & $1 $& $1 $\\
		\hline
		$\begin{bmatrix}
			1&\alpha&0 \\
			0 &1&\gamma \\
			0&0&1
		\end{bmatrix}$& $\epsilon^{\alpha(s-1)} $& $\epsilon^{\alpha(p-s+1)} $  & $\epsilon^{\{\alpha(s-1)+\gamma(t-1)\}} $& $\epsilon^{\{\alpha(p-s+1)+\gamma(p-t+1)\}} $\\
		\hline
		
		Dual irrep.	& \multicolumn{2}{c|}{$\rho_{2s-2}^*=\rho_{2s-1}$  $\& \, \, \rho_0=\rho_1$}& \multicolumn{2}{c|}{$\rho_{(t-1)p+2s-1}^*=\rho_{(t-1)p+2s}$} \\
		\hline	
		
		$\#$ irrep.	& \multicolumn{2}{c|}{$p$}& \multicolumn{2}{c|}{$p^2-p$} \\
		\hline
	\end{tabular} .
\end{table}

Since $1\leq  s,t \leq p$, and $\alpha, \gamma \in \mathbb{Z}_p$, so $\epsilon ^{(\alpha(s-1)+\gamma(t-1))}$ is generated by  $\epsilon$,  thus  the representation $\sigma_{(s-1,t-1)}$ maps an element  $g \in G$ to $\epsilon^{m},$ for some $m$, $1\leq  m \leq p.$

\subsection
{The non-abelian group $G_p$}
\begin{center}
	$G_p$ = $\left\{ (p\gamma,\delta)=\begin{bmatrix}
		1+p\gamma&\delta \\
		0 & 1 
	\end{bmatrix}\,|\,\, \gamma \in \mathbb{Z}_{p},\, \delta \in \mathbb{Z}_{p^{2}} \right\}$.
	\label{subs3.2}
\end{center}
Another presentation of this group is (\cite{Dummit}, pp 180)
\begin{center}
	$G_p= \langle  x,y\, | \, x^p=y^{p^2}=1 , xy=y^{p+1}x \rangle$.
\end{center}
The center of $G_p$ is 	$Z(G_p) =\left\langle y^p=\begin{bmatrix}
	1+p&0 \\
	0 & 1 
\end{bmatrix}\right\rangle$.
\begin{Note}
	For  the $p$ degree representations $\sigma_s$ of group $G_p$, if $y^p \in Z(G_p)$ then $\sigma_{2\eta-1}(y^p)= \rho_{p^2+2\eta-1}(y^p)=\omega^{\eta p^2}I_p$ and its dual $\sigma_{2\eta}(y^p)= \rho_{p^2+2\eta}(y^p)=\omega^{(p-\eta) p^2}I_p$,   $1\leq  \eta \leq  \frac{p-1}{2}$. Also 
	the elements  of $G_p - Z(G_p)$   map bijectively to the set
	\begin{center}
		$\sigma_s(G_p)-\sigma_s\big(Z(G_p)\big)$  = $\{A_{u} \in GL(p, \mathbb{F}) \,|\, \mbox{$Tr(A_{u}) = 0$ and $A_{u}^{p}= I_{p} $ or $A_{u}^{p^2}= I_{p} $ for $1 \leq u \leq p^3 -p$}\}$.
	\end{center}
	\label{note3.2}
\end{Note}

We have recorded the $p$ ordered irreducible representations of $G_p$ in the following table.
\begin{table}[h]
	\centering
	\caption{All irreducible representations of degree $p$ for the   group $G_p$.}
	
	\begin{tabular}{|c|c|c|c|c|c|c|c|}
			\hline
		All irrep. $\rightarrow$ &$\sigma_{2\eta-1}$&$\sigma_{2\eta}$\\
		\hline
		Runing variable $\rightarrow$	& \multicolumn{2}{c|}{1 $\leq \eta \leq \frac{p-1}{2}, $} \\
		\hline
		$g \in Z(G_p)$ & $\omega^{\eta p^2} I_p$ & $\omega^{(p-\eta)p^2} I_p$\\
		\hline
		
		$g \notin Z(G_p)$& See Note \ref{note3.2}&See Note \ref{note3.2}\\
		\hline
		
		Dual irrep.	& \multicolumn{2}{c|}{$\sigma_{2\eta-1}^*=\sigma_{2\eta}$} \\
		\hline
		$\#$ irrep.	& \multicolumn{2}{c|}{$p-1$} \\
		\hline
	\end{tabular}.
	\label{table3}
\end{table}

Thus the $p$ degree irreducible representations $\sigma_s$ is defined as below 	

$$ \sigma_{2 \eta-1}(x)	=\begin{bmatrix}
	0&1&0&\cdots&0\\
	0&0&1&\cdots&0\\
	0&0&0&\cdots&0\\
	\vdots& \vdots&\vdots&\ddots&\vdots\\
	1&0&0&\cdots&0
\end{bmatrix}^{\eta} \, \, and \, \sigma_{2 \eta-1}(y)	=\begin{bmatrix}
	\omega^{\eta p}&&0&\cdots&0\\
	0&\omega^{p^2+\eta p}&0&\cdots&0\\
	0&0&\omega^{2p^2+\eta p}&\cdots&0\\
	\vdots& \vdots&\vdots&\ddots&\vdots\\
	0&0&0&\cdots&\omega^{(p-1)p^2+\eta p}
\end{bmatrix},$$

$$ \sigma_{2 \eta}(x)	=\begin{bmatrix}
	0&1&0&\cdots&0\\
	0&0&1&\cdots&0\\
	0&0&0&\cdots&0\\
	\vdots& \vdots&\vdots&\ddots&\vdots\\
	1&0&0&\cdots&0
\end{bmatrix}^{-\eta} \,\,  and \,  \sigma_{2 \eta}(y)	=\begin{bmatrix}
	\omega^{p^2-\eta p}&&0&\cdots&0\\
	0&\omega^{2p^2-\eta p}&0&\cdots&0\\
	0&0&\omega^{3p^2-\eta p}&\cdots&0\\
	\vdots& \vdots&\vdots&\ddots&\vdots\\
	0&0&0&\cdots&\omega^{p^3-\eta p}
\end{bmatrix}.$$

\noindent There are $p^2$ representations $\sigma_{(s-1,t-1)},  1  \leq  s,t \leq p$ of degree  1.  As  $g  \in G_p - Z(G_p)$, we have $g=(p\gamma,\delta)=\begin{bmatrix}
	1+p\gamma&\delta \\
	0 & 1 
\end{bmatrix}$, $\delta \in \mathbb{Z}_{p^2}^*$,  $|G_p-Z(G_p)|=(p^2-1)p$. In the table 4 we have recorded all 1 degree representations of $G_p.$    \\

\begin{table}[h]
	\centering
	\caption{All irreducible representations of degree 1 for the   group $G_p$.}
	\begin{tabular}{|c|c|c|c|c|c|}
	\hline
	
	All irrep. $\rightarrow$&$\rho_{(t-1)p+2s-2}=$&$\rho_{(t-1)p+2s-1}=$&$\rho_{(t-1)p+2s-1}=$&$\rho_{(t-1)p+2s}=$\\
	& $\sigma_{(s-1,t-1)}$& $\sigma_{(p-s+1,p-t+1)}$ & $\sigma_{(s-1,t-1)}$& $\sigma_{(p-s+1,p-t+1)}$\\
	\hline
	
	Runing Variable	$\rightarrow$ &\multicolumn{2}{c|}{$1 \leq s \leq \frac{p+1}{2}, t=1$}& \multicolumn{2}{c|}{$1 \leq s \leq p,\,\, 2 \leq t \leq \frac{p+1}{2}$} \\
	\hline
		$y^{p \gamma}=\begin{bmatrix}
			1+p\gamma&0 \\
			0 & 1 
		\end{bmatrix}$& $\omega^{p^3}=1 $& $1 $ & $1 $& $1 $\\
		\hline
		$\begin{bmatrix}
			1+p\gamma&\delta \\
			0 & 1 
		\end{bmatrix}, \delta \neq 0$& $\omega^{p^2{\gamma(s-1)}} $& $\omega^{p^2\gamma(p-s+1)} $  & $\omega^{p^2\{\gamma(s-1)+\delta'(t-1)\}} $& $\omega^{p^2\{\gamma(p-s+1)+\delta'(p-t+1)\}} $\\
		&&& where $\delta' = \delta$ mod $p$& Where $\delta' = \delta$ mod $p$\\
		\hline
		
		Dual irrep.	& \multicolumn{2}{c|}{$\rho_{2s-2}^*=\rho_{2s-1}$  $\& \, \, \rho_0=\rho_1$}& \multicolumn{2}{c|}{$\rho_{(t-1)p+2s-1}^*=\rho_{(t-1)p+2s}$} \\
		\hline
		
		$\#$ irrep.	& \multicolumn{2}{c|}{$p$}& \multicolumn{2}{c|}{$p^2-p$} \\
		\hline
	\end{tabular} .
\end{table}
\pagebreak

Here  $\omega ^{p^2(\gamma(s-1)+\delta'(t-1))}$ is generated by $\omega^{p^2}$, thus  the representation $\sigma_{(s-1,t-1)}$ maps an element  $g \in G_p$ to $\omega^{mp^2}$,  for some $m$, $1 \leq  m \leq  p$.

\begin{Note}
	As a finite abelian group is finitely generated, here for the abelian groups $\mathbb{Z}_{p^3} $, $\mathbb{Z}_{p^2} \times \mathbb{Z}_{p} $ and $\mathbb{Z}_{p}  \times \mathbb{Z}_{p} \times \mathbb{Z}_{p} $ there exist the  finite generating subsets $\{a\} $, $\{b,c\} $ and $\{d,e,f\}  $ respectively. For example we may take $a=1+p^3\mathbb{Z}$, $b=(1+p^2\mathbb{Z}, 0+p\mathbb{Z})$, $c=(0+p^2\mathbb{Z},1+p\mathbb{Z})$, $d=(1+p\mathbb{Z},0+p\mathbb{Z},0+p\mathbb{Z})$, $e=(0+p\mathbb{Z},1+p\mathbb{Z},0+p\mathbb{Z})$ and $f=(0+p\mathbb{Z},0+p\mathbb{Z},1+p\mathbb{Z})$. Further order of  an element $g$ under $\rho_i, 1 \leq i \leq r$  is $|\rho_i(g)|$.
\end{Note}
\subsection{The cyclic group $\mathbb{Z}_{p^3}   $.}
Here $\mathbb{Z}_{p^3}   =\langle a \,|\,a^{p^3}=1  \rangle$.
The representation tables are given as below.

\begin{table}[h]
	\centering
	\caption{All irreducible representations of  the group ${\mathbb{Z}_{p^3}}$.}
	$T_{\mathbb{Z}_{p^3}} = $ 	
	\begin{tabular}{|c|c|c|c|c|c|c|c|c|}
		\hline
			All irrep. $\rightarrow$	& $\rho_t$&$\rho_{p^{s-1}+2t}$&$\rho_{p^{s-1}+2t+1}$\\
		\hline
		
		Runing variable $\rightarrow$	& &  \multicolumn{2}{c|}{$$ }\\	
		&$i=1$	 &  \multicolumn{2}{c|}{$  i \in \bigg\{x\,|\,\, gcd(p^3,xp^{3-s}) =p^{3-s}, \,1 \leq x \leq \frac{p(p^s-p^{s-1}-2)}{2(p-1)}+1$\bigg\}}\\
			&$t=1$ &  \multicolumn{2}{c|}{$t= i's$ place in above set }\\
		\hline
		$|\rho_i(a)|$	&1 & \multicolumn{2}{c|}{$p^s, \,s=1,2,3$}\\
		\hline

		$a$ & $\omega^{p^{3}}=1$&$\omega^{ip^{3-s}}, $&$\omega^{p^3-ip^{3-s}}$\\
		\hline
		
		Dual irrep.	&self &  \multicolumn{2}{c|}{$\rho_{p^{s-1}+2t}^*=\rho_{p^{s-1}+2t+1}$}\\
		\hline
		$\#$ irrep.	&1 &  \multicolumn{2}{c|}{$p^s-p^{s-1}$}\\
		\hline
	\end{tabular},
\end{table}

\begin{table}[h]
	\centering
	\caption{All irreducible representations of    the group ${\mathbb{Z}_{8}}$.}

	$T_{\mathbb{Z}_{8}} = $ 
	\begin{tabular}{|c|c|c|c|c|c|c|c|c|}
		\hline
		& $\rho_1$&$\rho_2$& $\rho_3$&$\rho_4$&$\rho_5$&$\rho_6$&$\rho_7$&$\rho_8$\\
		\hline
		$a$& $\omega^{8}$&$\omega^{4}$& $\omega^{2}$&$\omega^{6}$&$\omega$&$\omega^{7}$&$\omega^{3}$&$\omega^{5}$\\
		\hline
	\end{tabular} .
\end{table}
\subsection{ The group $\mathbb{Z}_{p^2} \times \mathbb{Z}_{p}  $. }
Here $\mathbb{Z}_{p^2} \times \mathbb{Z}_{p}  $ =$\langle a, b\,|\,a^{p^2}= b^{p}=1 ,ab= ba \rangle$.

\begin{table}[h!]
	\centering
	
	$T_{\mathbb{Z}_{p^2} \times \mathbb{Z}_{p}} = $ 
	
	\caption{All irreducible representations of   group ${\mathbb{Z}_{p^2} \times \mathbb{Z}_{p}}$.}
	\begin{tabular}{|c|c|c|c|c|c|c|c|c|c|c|c|c|}
		\hline
		All irrep.$\rightarrow$&$\sigma_{(s,t)}$	& $\sigma_{(p^2,t)}$&  $\sigma_{(p^2,p-t)}$& $\sigma_{(s,p)}$&  $\sigma_{(p-s,p)}$ & $\sigma_{(s,t)}$&  $\sigma_{(p-s,p-t)}$
		& $\sigma_{(s,t)}$&  $\sigma_{(p^2-s,p-t)}$ & $\sigma_{(s,p)}$&  $\sigma_{(p-s,p)}$\\
		\hline
		Runing &$s$ & \multicolumn{2}{c|}{$$} & %
		\multicolumn{2}{c|}{$$} & \multicolumn{2}{c|}{$ t \in \{1,..,p-1\}$}&
		\multicolumn{2}{c|}{$ t \in \{1,..,p-1\}$} & \multicolumn{2}{c|}{$$}\\
		variable & $=t$& \multicolumn{2}{c|}{$$} & %
		\multicolumn{2}{c|}{$$} & \multicolumn{2}{c|}{$ \& \, 1 \leq s\leq \frac{p-1}{2} \, or $}&
		\multicolumn{2}{c|}{$ \& \, 1 \leq s\leq \frac{p^2-p}{2}\, or$} & \multicolumn{2}{c|}{$$}\\
		
		$\rightarrow$	&$=p$ & \multicolumn{2}{c|}{$1 \leq t \leq \frac{p-1}{2}$} & %
		\multicolumn{2}{c|}{$1 \leq s \leq \frac{p-1}{2}$} & \multicolumn{2}{c|}{$ s\in\{1,...,p-1\}$}&
		\multicolumn{2}{c|}{$  s \in \{1,...,p^2-p\}$} & \multicolumn{2}{c|}{$1 \leq s \leq \frac{p^2-p}{2}$}\\
	& & \multicolumn{2}{c|}{$$} & %
		\multicolumn{2}{c|}{$$} & \multicolumn{2}{c|}{$ \& \, 1\leq t \leq \frac{p-1}{2}$}&
		\multicolumn{2}{c|}{$ \& \,1 \leq t \leq \frac{p-1}{2}$} & \multicolumn{2}{c|}{$$}\\
		\hline
		$|\rho_i(a)|$&1 & \multicolumn{2}{c|}{$1$} & %
		\multicolumn{2}{c|}{$p$} & \multicolumn{2}{c|}{$p$}&
		\multicolumn{2}{c|}{$p^2$} & \multicolumn{2}{c|}{$p^2$}\\
		
		\hline
		$|\rho_i(b)|$&1 & \multicolumn{2}{c|}{$p$} & %
		\multicolumn{2}{c|}{$1$} & \multicolumn{2}{c|}{$p$}&
		\multicolumn{2}{c|}{$p$} & \multicolumn{2}{c|}{$1$}\\

		\hline
		$a$ &1& 1&1& $\omega^{sp^2}$&$\omega^{(p-s)p^2}$ & $\omega^{sp^2}$&$\omega^{(p-s)p^2}$& $\omega^{sp}$&$\omega^{(p^2-s)p}$ & $\omega^{sp}$&$\omega^{(p^2-s)p}$\\
	\hline
		$b$ &1&  $\omega^{tp^{2}}$&$\omega^{(p-t)p^{2}}$&1 &1& $\omega^{sp^2}$&$\omega^{(p-t)p^2}$ & $\omega^{tp^2}$&$\omega^{(p-t)p^2}$ & 1&1\\
		\hline

		Dual irrep.	&self & \multicolumn{2}{c|}{$\sigma_{(p^2,t)}^*$  $  = \sigma_{(p^2,p-t)}$} & %
		\multicolumn{2}{c|}{ $\sigma_{(s,p)}^*$ $=\sigma_{(p-s,p)}$} & \multicolumn{2}{c|}{$\sigma_{(s,t)}^* =$ $\sigma_{(p-s,p-t)}$}&
		\multicolumn{2}{c|}{$\sigma_{(s,t)}^*=$ $\sigma_{(p^2-s,p-t)}$} & \multicolumn{2}{c|}{$\sigma_{(s,p)}^*=$ $\sigma_{(p-s,p)}$}\\
		\hline
		
		$\#$ irrep.	&1 & \multicolumn{2}{c|}{$p -1$} & %
		\multicolumn{2}{c|}{$p-1$} & \multicolumn{2}{c|}{$(p-1)(p-1)$}&
		\multicolumn{2}{c|}{$(p^2-p)(p-1)$} & \multicolumn{2}{c|}{$p^2-p$}\\
		\hline
	\end{tabular},
\end{table}

\begin{table}[h!]
	\centering
	\caption{All irreducible representations of   group ${\mathbb{Z}_{4} \times \mathbb{Z}_{2}}$.}

	$T_{\mathbb{Z}_{4} \times \mathbb{Z}_{2}} = $ 
	\begin{tabular}{|c|c|c|c|c|c|c|c|c|}
		\hline
		& $\rho_1$&$\rho_2$& $\rho_3$&$\rho_4$&$\rho_5$&$\rho_6$&$\rho_7$&$\rho_8$\\
		\hline
		$a$& $\omega^{8}$&$\omega^{8}$& $\omega^{4}$&$\omega^{4}$&$\omega^{2}$&$\omega^{6}$&$\omega^{2}$&$\omega^{6}$\\
		\hline
		$b$& $\omega^{8}$&$\omega^{4}$& $\omega^{8}$&$\omega^{4}$&$\omega^{4}$&$\omega^{4}$&$\omega^{8}$&$\omega^{8}$\\
		\hline
	\end{tabular} .
	
\end{table}
\subsection{ The group $\mathbb{Z}_{p} \times \mathbb{Z}_{p}  \times \mathbb{Z}_{p}. $ }
Here $\mathbb{Z}_{p} \times \mathbb{Z}_{p}  \times \mathbb{Z}_{p}   =\langle a, b,c\,|\,a^{p}= b^{p}=c^{p}=1 ,ab= ba,ac=ca,bc=cb \rangle$. The corresponding tables are given by

\begin{table}[h!]
	\centering
	$T_{\mathbb{Z}_{p} \times \mathbb{Z}_{p} \times \mathbb{Z}_{p}} = $
	\caption{All irreducible representations of $\mathbb{Z}_{p} \times \mathbb{Z}_{p} \times \mathbb{Z}_{p}$ } 
	\begin{tabular}{|c|c|c|c|c|c|c|c|}
		\hline
			3-tuples $\rightarrow$ 	& \multicolumn{2}{c|}{ All place are same } & %
		\multicolumn{2}{c|}{Exactly two place are same } & \multicolumn{2}{c|}{ All place are distinct}\\
	 $ (s,t,m)$  	& \multicolumn{2}{c|}{ $(s,s,s) $} & %
		\multicolumn{2}{c|}{$(s,s,m), (s,m,s),(m,s,s), m \neq s$ } & \multicolumn{2}{c|}{ $(s,t,m)$, $s \neq t \neq m \neq s$}\\
		\hline
		All irrep. $\rightarrow$& $\sigma_{(s,s,s)}$& $\sigma_{(p-s,p-s,p-s)}$& $\sigma_{(s,s,m)}$&$\sigma_{(p-s,p-s,p-m)}$ & $\sigma_{(s,t,m)}$&$\sigma_{(p-s,p-t,p-m)}$\\
		\hline
		Runing& \multicolumn{2}{c|}{$\frac{p+1}{2} \leq s \leq p$ } & %
		\multicolumn{2}{c|}{$s=p \,\, \& \, 1 \leq m \leq \frac{p-1}{2}$ or } & \multicolumn{2}{c|}{ 	$1 \leq s,t,m \leq p$, }\\
	variable $\rightarrow$	& \multicolumn{2}{c|}{} & %
		\multicolumn{2}{c|}{$1 \leq s \leq \frac{p-1}{2}\,  \& \, m \in \{1,2,\cdots p\}$} & \multicolumn{2}{c|}{ 	$$}\\

		\hline
		$a$& $\omega^{sp^{2}}$ & $\omega^{(p-s)p^{2}}$& $\omega^{sp^{2}}$ & $\omega^{(p-s)p^{2}}$ & $\omega^{sp^{2}}$ & $\omega^{(p-s)p^{2}}$\\
		\hline
		
		$b$ & $\omega^{sp^{2}}$ & $\omega^{(p-s)p^{2}}$& $\omega^{sp^{2}}$ & $\omega^{(p-s)p^{2}}$ & $\omega^{tp^{2}}$ & $\omega^{(p-t)p^{2}}$\\
		\hline
		
		$c$ & $\omega^{sp^{2}}$ & $\omega^{(p-s)p^{2}}$& $\omega^{mp^{2}}$ & $\omega^{(p-m)p^{2}}$ & $\omega^{mp^{2}}$ & $\omega^{(p-m)p^{2}}$\\
		\hline
		Twice & \multicolumn{2}{c|}{$\sigma_{(p,p,p)}=\sigma_{(0,0,0)} $} & %
		\multicolumn{2}{c|}{} & \multicolumn{2}{c|}{$\sigma_{(s,t,m)}= \sigma_{(p-s_1,p-t_1,p-m_1) } \, if $}\\
			repeated irrep. & \multicolumn{2}{c|}{Only trivial irrep.} & %
		\multicolumn{2}{c|}{} & \multicolumn{2}{c|}{$s+s_1=t+t_1=m+m_1=p $\, or}\\
	 $\rightarrow$	& \multicolumn{2}{c|}{repeated twice} & %
		\multicolumn{2}{c|}{No  irrep. repeated twice } & \multicolumn{2}{c|}{ $2p$, Each irrep. repeated twice}\\
		\hline
		
		Dual irrep. 	& \multicolumn{2}{c|}{$\sigma_{(s,s,s)}^*=\sigma_{(p-s,p-s,p-s)}$} & %
		\multicolumn{2}{c|}{$\sigma_{(s,s,m)}^*= \sigma_{(p-s,p-s,p-m)}$} & \multicolumn{2}{c|}{$\sigma_{(s,t,m)}^*= \sigma_{(p-s,p-t,p-m)}$}\\
		\hline
		Distinct	& \multicolumn{2}{c|}{$$} & %
		\multicolumn{2}{c|}{$$} & \multicolumn{2}{c|}{$$}\\
		
		$\#$ irrep.	& \multicolumn{2}{c|}{$p$} & %
		\multicolumn{2}{c|}{$3p(p-1)$} & \multicolumn{2}{c|}{$p(p-1)(p-2)$}\\
		\hline
	\end{tabular},
\end{table}

\begin{table}[h!]
	\centering
	
	\caption{All irreducible representations of $\mathbb{Z}_{2} \times \mathbb{Z}_{2} \times \mathbb{Z}_{2}$}
	$T_{\mathbb{Z}_{2} \times \mathbb{Z}_{2} \times \mathbb{Z}_{2}} = $
	\begin{tabular}{|c|c|c|c|c|c|c|c|c|}
		\hline
		& $\rho_1$&$\rho_2$& $\rho_3$&$\rho_4$&$\rho_5$&$\rho_6$&$\rho_7$&$\rho_8$\\
		\hline
		$a$& $\omega^{8}$&$\omega^{4}$& $\omega^{8}$&$\omega^{8}$&$\omega^{8}$&$\omega^{4}$&$\omega^{4}$&$\omega^{4}$\\
		\hline
		
		$b$& $\omega^{8}$&$\omega^{4}$& $\omega^{4}$&$\omega^{4}$&$\omega^{8}$&$\omega^{4}$&$\omega^{8}$&$\omega^{8}$\\
		\hline
		
		$c$& $\omega^{8}$&$\omega^{4}$& $\omega^{8}$&$\omega^{4}$&$\omega^{4}$&$\omega^{8}$&$\omega^{8}$&$\omega^{4}$\\
		\hline
		
	\end{tabular}.
	
\end{table}
\pagebreak
\begin{Note}
	An irreducible representation   $\rho_{2i}$ is seated together (preceded or succeeded) to its dual  $\rho_{2i+1}$     in $\rho$. 
	\label{note3.3}
\end{Note}

\noindent Now 
\begin{equation}
	\rho= k_{1} \rho_{1}\oplus k_{2} \rho_{2}\oplus  ............. \oplus k_{r} \rho_{r},
	\label{directsum}
\end{equation}
where for every $1 \leq i \leq r$,  $k_{i} \rho_{i}$ stands for the direct sum of  $k_{i}$ copies of the irreducible representation $\rho_{i}$.

\noindent Let $\chi$ be the corresponding character of the representation $\rho$, then  
$$ \chi= k_{1} \chi_{1}+ k_{2} \chi_{2}+  ............. + k_{r} \chi_{r},$$
where $\chi_{i}$ is the character of $\rho_{i}$, $\forall$ $1 \leq i \leq r$.
Dimension of the character $\chi$ is being calculated at the identity element of the group. i.e,
\[ dim(\rho) = \chi(1)=tr(\rho(1)).\]
\begin{equation}
	\implies d_{1} k_{1} + d_{2}k_{2} +............. +d_{r}k_{r}= n. 
	\label{splitn}
\end{equation}
\begin{Note}
	Equation (\ref{splitn}) holds in more general situation, which  helps us to find all  possible distinct r-tuples ($k_{1}, k_{2}, ......, k_{r}$), which correspond to the distinct n degree  representations (up to isomorphism) of a given finite group.
\end{Note}

\begin{theorem}
	Let $G$ be a group  of order $p^3$ with $p$ a prime. If $\sigma$ is an irreducible representation of $G$ of degree $p$, then $\sigma$ is a faithful representation.
\end{theorem}
\begin{proof}
	For  non-abelian groups it is clear from the tables $\ref{table1}$ and $\ref{table3}$  in the subsections   $\ref{subs3.1}$ and $\ref{subs3.2}$	   in this paper if $p$ is an odd prime, whereas for $p=2$ it follows from  the tables  $T_{D_{4}}$ and $T_{Q_{8}}$  in the  article \cite{DCT}. For an abelian group there is no irreducible representation of degree $p$.
\end{proof}
\section{Existence of non-degenerate invariant   forms.}
\indent An element in  the space of invariant  bilinear forms under representation of  a finite  group is either    non-degenerate or  degenerate. All the elements of the space are degenerate when  $k_{2i} \neq k_{2i+1}$,  such a space is called a degenerate invariant space which has  also been discussed in \cite{DCT} for the groups of order $8$. How many such representations exist out of  total representations, is a matter of investigation. Some of the spaces contain both non-degenerate and degenerate invariant  bilinear forms under   a particular representation. In this section we compute the number of such representations of the group $G$ of order $p^3$, with $p$ an odd prime.\\
\begin{Remark}
	The space $\Xi_{G}'$  of invariant bilinear forms under an n degree representation $\rho$ contains only those $X \in \mathbb{M}_{n}(\mathbb{F})$ whose  $(i,j)^{th}$ block is a $\text{\huge0}$ sub-matrix of order ${d_{i}k_{i}\times d_{j}k_{j}}$ when $(i,j) \notin A_G =\{(i,j)\,\,|\,\,$ \mbox{ $\rho_i$ and $\rho_j$ are dual to each other} \}  whereas  for  $(i,j) \in A_G$ with $d_i=d_j=1$, the block matrix $X^{ij}_{d_ik_i \times d_jk_j}$ is given by
	$$X^{ij}_{d_ik_i \times d_jk_j}= X^{ij}_{k_i \times k_j}=
	\begin{bmatrix}
		x^{ij}_{1,1}&x^{ij}_{1,2}&\cdots&x^{ij}_{1, k_{j}}\\
		x^{ij}_{2,1}&x^{ij}_{2,2}&\cdots&x^{ij}_{2,k_{j}}\\
		\vdots & \vdots &\ddots&\vdots\\
		
		x^{ij}_{k_{i},1}&x^{ij}_{k_{i},2}&\cdots&x^{ij}_{k_{i},k_{j}}
		\label{remark4.1}
	\end{bmatrix}.$$


	And for  $(i,j) \in A_{G_{p}}$ with $d_i=d_j=p$ it is, 
	$$X^{ij}_{d_ik_i \times d_jk_j}= X^{ij}_{pk_i \times pk_j}=
	\begin{bmatrix}
		x^{ij}_{1,p}L&x^{ij}_{1,2p}L&\cdots&x^{ij}_{1,pk_j}L\\
		x^{ij}_{p+1,p}L&x^{ij}_{p+1,2p}L&\cdots&x^{ij}_{p+1,pk_j}L\\
		\vdots & \vdots &\ddots&\vdots\\
		
		x^{ij}_{(k_i-1)p+1,p}L&	x^{ij}_{(k_i-1)p+1,2p}L&\cdots&	x^{ij}_{(k_i-1)p+1,pk_j}L
		
	\end{bmatrix},$$
	$$ where \, L=\begin{bmatrix}
		0&\cdots&0&1\\
		0&\cdots&1&0\\
		\vdots & \reflectbox{$\ddots$} &\vdots&\vdots\\
		
		1&\cdots&0&0
		
	\end{bmatrix}_{p \times p}.$$
	
	\label{remark}
\end{Remark} 
\begin{lemma}
	If $X \in \Xi_{G}'$,  and $k_{2i} \neq k_{2i+1}$, then $X$ must be singular.
	\label{singular}
\end{lemma}  

\begin{proof}  With reference to the above remark and Note \ref{note3.3}, for every $X \in \Xi_{G}'$, we have $X= [X^{ij}_{d_ik_i \times d_jk_j}]_{(i,j) \in A_G}$. I.e.
	
	$$X= 
	\begin{bmatrix}
		X^{11}_{d_1k_1 \times d_1k_1}& $\mbox{\Huge 0} $ & \cdots &$\mbox{\Huge 0}$\\
		$\mbox{\Huge 0}$&\begin{bmatrix}
			$\text{\huge0}$&X^{23}_{d_2k_{2} \times d_3k_3} \\
			X^{32}_{d_3k_{3}\times d_2 k_2}&$\text{\huge0}$
		\end{bmatrix}& \cdots &$\mbox{\Huge 0}$\\
		\vdots & \vdots &\ddots&\vdots\\
		
		$\mbox{\Huge 0}$&$\mbox{\Huge 0}$&\cdots& \begin{bmatrix}
			$\text{\huge0}$&X^{r-1\, \, r}_{d_{r-1}k_{r-1}\times d_rk_r} \\
			X^{r\,\,r-1}_{d_rk_{r}\times d_{r-1}k_{r-1}}&$\text{\huge0}$
		\end{bmatrix} \\
	\end{bmatrix},$$

	with $X^{ij}_{d_{i}k_{i} \times d_jk_j}=  C_{k_{i}\rho_{i}(g)}^{t}X^{ij}_{d_ik_i \times d_{j}k_{j}}C_{k_{i}\rho_{i}(g)}$ , for $(i,j) \in A_G$. If  $k_{2i}\neq k_{2i+1}$ and  since   $d_{2i}=d_{2i+1}$, so  the number of  rows  and columns of $X^{ij}_{d_{i}k_{i} \times d_jk_j}$ is differ  hence  either rows (or columns) is linearly dependent.   Thus the result follows. 
\end{proof} 

In the next lemma we characterise the representations of $G$ each of which admits a non-degenerate invariant bilinear form.
To prove  the next lemma we will choose  only those $X \in \mathbb{M}_{n}(\mathbb{F})$ whose $(i,j)^{th} $ block is zero for $(i,j)\notin A_G$, whereas  for  $(i,j) \in A_G$ with $k_i=k_j$ and  the block matrices  $X^{ij}_{d_{i}k_{i} \times d_jk_j}$, is   non-singular. 

\begin{lemma} For $n \in \mathbb{Z}^{+}$, an $n$-degree representation of $G$ has a non-degenerate invariant  bilinear form iff $k_{2i}=k_{2i+1}, 1 \leq i \leq \frac{r-1}{2}$.
	\label{nondeg}
\end{lemma}
\begin{proof} From equation (\ref{splitn}) we have $d_{1} k_{1} + d_{2}k_{2} +............. +d_{r}k_{r}= n$ and  chose $X \in \mathbb{M}_{n}(\mathbb{F})$ such that 			$$X= 
	\begin{bmatrix}
		X^{11}_{d_1k_1 \times d_1k_1}& $\mbox{\Huge 0} $ & \cdots &$\mbox{\Huge 0}$\\
		$\mbox{\Huge 0}$&\begin{bmatrix}
			$\text{\huge0}$&X^{23}_{d_2k_{2} \times d_3k_3} \\
			X^{32}_{d_3k_{3}\times d_2 k_2}&$\text{\huge0}$
		\end{bmatrix}& \cdots &$\mbox{\Huge 0}$\\
		\vdots & \vdots &\ddots&\vdots\\
		
		$\mbox{\Huge 0}$&$\mbox{\Huge 0}$&\cdots& \begin{bmatrix}
			$\text{\huge0}$&X^{r-1\, \, r}_{d_{r-1}k_{r-1}\times d_rk_r} \\
			X^{r\,\,r-1}_{d_rk_{r}\times d_{r-1}k_{r-1}}&$\text{\huge0}$
		\end{bmatrix} \\
	\end{bmatrix}.$$
	Suppose $k_{2i}=k_{2i+1}$, and for $(i,j) \in A_G$, the chosen  block matrices 
	$X^{ij}_{d_{i}k_{i} \times d_{j}k_{j}}$ is non-singular. Thus the rows (or columns) of the $X$ is linearly independent. Also from remark \ref{remark4.1} we have  $X^{ij}_{d_{i}k_{i} \times d_jk_j}=  C_{k_{i}\rho_{i}(g)}^{t}X^{ij}_{d_ik_i \times d_{j}k_{j}}C_{k_{i}\rho_{i}(g)}$ so  $C_{\rho(g)} ^{t}X C_{\rho(g)}= X$. Therefore  $X \in \Xi_{G}'$.\\
	On other hand $X \in \Xi_{G}'$ non-singular implies that $X^{ij}_{d_{i}k_{i} \times d_{j}k_{j}}$ is non-singular for $(i,j) \in A_G$, therefore the corresponding block matrix is square which is possible only when  $k_{2i}=k_{2i+1},$ $\forall \,i \in \{1,2,...,\frac{r-1}{2}\}$. 
\end{proof}
Note that the Lemmas \ref{singular} and \ref{nondeg} can be covered in a more general situation by stating as 'no non-trivial irreducible representation of a finite $p$-group can be self dual'. For if $L$ a finite $p$-group and $V$ a non-trivial irreducible representation of $L$, replacing $L$ by its image in the matrix group
(the general linear group), we may assume that the representation is faithful. Being a non-trivial finite $p$-group, $L$ has non-trivial centre. Let $g$ be a non-trivial central element in $L$. The action of $g$ on $V$ is by multiplication by a root of unity $\zeta\neq \pm  1 $. Its action on the dual of $V$ is by multiplication by $\bar\zeta (=\zeta^{-1})$. Since $\zeta\not=\bar\zeta$, it follows that $V$ is not self-dual.

\begin{corollary} For $n \in \mathbb{Z}^{+}$, an n degree representation of a group of order $p^3$ has a non-degenerate invariant bilinear form if and only if every irreducible representation and its dual have same multiplicity in the representation. 
	\label{lemma43}
\end{corollary}
\begin{proof} Follows from the proof of the lemma $\ref{nondeg}$ .
\end{proof}


\begin{Remark} If $\mathbb{F}$ is algebraically closed, it has infinitely many non zero elements, hence  if there is one non-degenerate  invariant bilinear form in the space $\Xi_{G}$, it has   infinitely many.\\
\end{Remark}

\begin{lemma}
	Let $G$ be a group of order $p^3$ and $n\in\mathbb{N}$. Then the  number of  $n$ degree representations of $G$ each of which admits a non-degenerate invariant  bilinear form is  $\displaystyle\sum_{\ell=0}^{ \lfloor \frac{n}{2p} \rfloor}\Bigg{[}\binom{\ell+\frac{p-3}{2}}{\frac{p-3}{2}} \sum_{s=0}^{\lfloor \frac{n-2p\ell}{2} \rfloor} \binom{s+\frac{p^2-3}{2}}{\frac{p^2-3}{2}}\Bigg]$ when $G$ is non-abelian whereas it is  $\displaystyle \sum_{\ell=0}^{\lfloor \frac{n}{2}\rfloor}\binom{\ell+\frac{p^{3}-3}{2}}{\frac{p^{3}-3}{2}}$ when $G$ is abelian.
	\label{lemma3.9}
\end{lemma}

\begin{proof} 
	
	Let $\rho = \oplus_{i=1}^{r} k_{i}\rho_{i}$ be an $n$ degree representation of $G$ which admits a non-degenerate bilinear form. So, we have  $k_{2i}=k_{2i+1}$, $1 \leq i \leq \frac{r-1}{2}$. In the non-abelian case,  $G$ is either $G_{p}$ or $Heis(\mathbb{Z}_{p})$ and for each of these two, we have  $r = p^{2}+p-1$,  $d_{i}$ =1 for $1\leq i \leq p^{2}$ and $d_{i}=p$ for $p^{2}+1 \leq i \leq p^{2}+p-1$. Now from equation (\ref{splitn}), we have
	
	\[
	k_1+2( \underbrace{k_2+k_4+ \cdots + k_{p^2-1}}_{\frac{p^2-1}{2}} ) +2p( \underbrace{k_{p^2+1}+k_{p^2+3}+ \cdots + k_{p^2+p-2}}_{\frac{p-1}{2}} ) =n
	\] 
	\[\implies
	k_1+2( \underbrace{k_2+k_4+ \cdots + k_{p^2-1}}_{\frac{p^2-1}{2}} ) =n -2 p( \underbrace{k_{p^2+1}+k_{p^2+3}+ \cdots + k_{p^2+p-2}}_{\frac{p-1}{2}} )
	\] 
	\begin{equation}
		\implies	k_1+2( \underbrace{k_2+k_4+ \cdots + k_{p^2-1}}_{\frac{p^2-1}{2}} ) =n -2p\ell.
		\label{equation3}
	\end{equation}

	To solve the above  equation we have $( \underbrace{k_{p^2+1}+k_{p^2+3}+ \cdots + k_{p^2+p-2}}_{\frac{p-1}{2}} )=\ell$, $0\leq \ell \leq \lfloor \frac{n}{2p}\rfloor$. i.e,  we have $\lfloor \frac{n}{2p}\rfloor +1$ equations. The $\ell^{th}$ equation is    \\
	
	\begin{equation}
		\underbrace{k_{p^2+1}+k_{p^2+3}+ \cdots + k_{p^2+p-2}}_{\frac{p-1}{2}}=\ell.
		\label{eqn4}
	\end{equation}

	The number of distinct solution to  above equations \ref{eqn4}   is $\binom{\ell+\frac{p-1}{2}-1}{\frac{p-1}{2}-1}$, $0 \leq \ell \leq \lfloor \frac{n}{2p}\rfloor$.\\
	Further	from equation \ref{equation3} we have
	\[
	k_1= n -2p\ell-2( \underbrace{k_2+k_4+ \cdots + k_{p^2-1}}_{\frac{p^2-1}{2}} )
	\] 
	\[\implies
	k_1= n -2p\ell-2\lambda,
	\] 
	where   $\underbrace{k_2+k_4+ \cdots + k_{p^2-1}}_{\frac{p^2-1}{2}} = \lambda $, $0 \leq \lambda \leq \lfloor \frac{n-2p\ell}{2} \rfloor $. i.e, we have $\lfloor \frac{n-2p\ell}{2} \rfloor+1$ equations and the number of solutions for every $ \lambda$ to the  equation is $\binom{\lambda +\frac{p^2-1}{2}-1}{\frac{p^2-1}{2}-1}.$\\
	Thus the number of all  distinct    $p^{2}+p-1$ tuples ($k_{1}, k_{2}, k_{3} ......., k_{p^{2}+p-2}, k_{p^{2}+p-1}$)  with $k_{2i} =k_{2i+1}$ is\\ $$\sum_{\ell=0}^{ \lfloor \frac{n}{2p} \rfloor}\Bigg{[}\binom{\ell+\frac{p-1}{2}-1}{\frac{p-1}{2}-1} \sum_{\lambda=0}^{\lfloor \frac{n-2p\ell}{2} \rfloor} \binom{\lambda+\frac{p^2-1}{2}-1}{\frac{p^2-1}{2}-1} \Bigg].$$\\
	\vskip1mm
	Now in the abelian case $G$ is either of $\mathbb{Z}_{p} \times \mathbb{Z}_{p}\times \mathbb{Z}_{p}$,  $\mathbb{Z}_{p^2} \times \mathbb{Z}_{p}$ and $\mathbb{Z}_{p^3} $  for each of which  $r=p^{3}$ and  $d_{i} =1$ for 
	$1 \leq i \leq p^{3}$. Now from  (\ref{splitn}), we have 
	$$ k_{1} + 2(\underbrace{ k_{2} +k_4+.....................+k_{p^{3}-1}}_{\frac{p^3-1}{2}})= n. $$
	Thus  the number of all  distinct  $p^3$-tuples ($k_{1}, k_{2}, k_{3},\cdots \cdots \cdots , k_{p^3}$) with $k_{2i}=k_{2i+1}$ is $\sum_{\ell=0}^{\lfloor \frac{n}{2}\rfloor}\binom{\ell+\frac{p^{3}-1}{2}-1}{\frac{p^{3}-1}{2}-1}$.\\

	Thus from equation (\ref{splitn}) and Theorem $\ref{theorem2.2}$ the  number of  $n$ degree representations (upto isomorphism)  of a group $G$ consisting non-degenerate invariant bilinear form is $\sum_{\ell=0}^{ \lfloor \frac{n}{2p} \rfloor}\Bigg{[}\binom{\ell+\frac{p-1}{2}-1}{\frac{p-1}{2}-1} \sum_{\lambda =0}^{\lfloor \frac{n-2p\ell}{2} \rfloor} \binom{\lambda+\frac{p^2-1}{2}-1}{\frac{p^2-1}{2}-1} \Bigg]$ for non-abelian groups  and  $\sum_{\ell=0}^{\lfloor \frac{n}{2}\rfloor}\binom{l+\frac{p^{3}-1}{2}-1}{\frac{p^{3}-1}{2}-1}$ for abelian groups of order $p^3$, with $p$ an odd prime.

\end{proof}

\begin{definition}
	The space $\Xi_{G}$ of invariant bilinear forms  is   called  degenerate  if it's all elements are degenerate.
\end{definition}
We  will  discuss about the degenerate  invariant  space in the later section.\\

\section{Dimensions of  spaces of  invariant bilinear forms under the representations of  groups of order $p^3$ with prime $p>2$.}
\noindent  The  space of invariant   bilinear forms under an n degree representation is generated by finitely many vectors, so its  dimension is finite along with its  symmetric  and the skew-symmetric subspace. In this section  we formulate the dimension of the space of invariant bilinear forms under a representation of a group of order $p^3$, with $p$ an odd prime.\\
\begin{theorem}
	If $\Xi_{G}$ is the  space of invariant bilinear forms  under an n degree representation  $\rho =\oplus_{i=1}^{r}k_{i}\rho_{i}$ of a group $G$ of order $p^3$, then  dim$(\Xi_{G})= \sum_{(i,j) \in A_G}k_{i}k_{j}$. 	 
	\label{th5.1}
\end{theorem}
\begin{proof}
	For every $X \in \Xi_{G}'$, we have 
	
	$$X= 
	\begin{bmatrix}
		X^{11}_{d_1k_1 \times d_1k_1}& $\mbox{\Huge 0} $ & \cdots &$\mbox{\Huge 0}$\\
		$\mbox{\Huge 0}$&\begin{bmatrix}
			$\text{\huge0}$&X^{23}_{d_2k_{2} \times d_3k_3} \\
			X^{32}_{d_3k_{3}\times d_2 k_2}&$\text{\huge0}$
		\end{bmatrix}& \cdots &$\mbox{\Huge 0}$\\
		\vdots & \vdots &\ddots&\vdots\\
		
		$\mbox{\Huge 0}$&$\mbox{\Huge 0}$&\cdots& \begin{bmatrix}
			$\text{\huge0}$&X^{r-1\, \, r}_{d_{r-1}k_{r-1}\times d_rk_r} \\
			X^{r\,\,r-1}_{d_rk_{r}\times d_{r-1}k_{r-1}}&$\text{\huge0}$
		\end{bmatrix} \\
	\end{bmatrix},$$

	with $X^{ij}_{d_{i}k_{i} \times d_jk_j}=  C_{k_{i}\rho_{i}(g)}^{t}X^{ij}_{d_{i}k_{i} \times d_jk_j}C_{k_{j}\rho_{j}(g)}$, for $(i,j) \in A_G$ and to generate each of these blocks of $X$ it needs $k_ik_j$ vectors from $\Xi_G'$. Thus  the result follows.	
\end{proof}

\begin{corollary}
	The  space of invariant symmetric bilinear forms  under an n degree representation  $\rho =\oplus_{i=1}^{r}k_{i}\rho_{i}$ of a group $G$ of order $p^3$ has dimension $=$ $\frac{k_{1}(k_{1}+1)}{2} +\sum_{\substack{(i,j) \in {A_G} \\ i \neq j}}\frac{k_{i}k_{j}}{2}$.
\end{corollary}
\begin{proof}
	Follows  from the proof of theorem $\ref{th5.1}$ .
\end{proof}

\begin{corollary}
	The space of invariant skew-symmetric bilinear forms  under an n degree representation  $\rho =\oplus_{i=1}^{r}k_{i}\rho_{i}$ of a group $G$ of order $p^3$ has  dimension $=$ $\frac{k_{1}(k_1-1)}{2} +\sum_{\substack{(i,j) \in {A_G} \\ i \neq j}}\frac{k_{i}k_{j}}{2}$.
\end{corollary}

\begin{proof}
	Follows  from the proof of theorem $\ref{th5.1}$ .
\end{proof}


\section{Main results \label{proof1.1}}
Here we present the proofs of the  main theorems stated in the Introduction section.
\begin{flushleft}
	\textbf{Proof of Theorem $\ref{theorem1.1}$}
	Since $G$ is the group of order $p^{3}$,  with an odd prime $p$ and degree of the representation $\rho$	is $n$,  if $G$ is either $G_{p}$ or $Heis(\mathbb{Z}_{p})$, we have  $r = p^{2}+p-1$,  $d_{i}$ =1 for $1\leq i \leq p^{2}$ and $d_{i}=p$ for $p^{2}+1 \leq i \leq p^{2}+p-1$. Now from equation (\ref{splitn}), we have 
	$$ \hspace{1.2cm} k_{1} + k_{2} +..........+k_{p^{2}} +pk_{p^{2}+1} +........+pk_{p^{2}+p-1}= n $$
	$$ Or, \hspace{1cm} k_{1} + k_{2} +..........+k_{p^{2}} = n- p(k_{p^{2}+1} +........+k_{p^{2}+p-1}). $$
	\begin{equation}
		Or, \hspace{1cm} k_{1} + k_{2} +..........+k_{p^{2}} = n- p\mu,
	\end{equation}
	
	where $\mu=k_{p^{2}+1} +........+k_{p^{2}+p-1}$, $ 0 \leq \mu \leq \lfloor \frac{n}{p} \rfloor$. i.e,  we have $\lfloor \frac{n}{p}\rfloor +1$ equations placed in the chronological order and the $\mu^{th}$ equation is given by\\
	
	\begin{equation}
		k_{p^{2}+1} + k_{p^{2}+2}+........+k_{p^{2}+p-1}= \mu.
		\label{equation6}
	\end{equation}

	The number of distinct solutions to   equation \ref{equation6} is $\binom{\mu+p-2}{p-2}$, $0 \leq \mu \leq \lfloor \frac{n}{p} \rfloor$.\\
	Thus the number of all  distinct    $p^{2}+p-1$ tuples ($k_{1}, k_{2}, ......., k_{p^{2}+p-1}$) is $\sum_{\mu=0}^{ \lfloor \frac{n}{p}\rfloor} \binom{\mu +p-2}{p-2}\binom{n-\mu p+p^{2}-1}{p^{2}-1}$.\\
	\vskip1mm
	On the otherhand if $G$ is either of $\mathbb{Z}_{p} \times \mathbb{Z}_{p}\times \mathbb{Z}_{p}$,  $\mathbb{Z}_{p^2} \times \mathbb{Z}_{p}$ and $\mathbb{Z}_{p^3} $  then  $r=p^{3}$ and  $d_{i}$ =1 for 
	$1 \leq i \leq p^{3}$. Now from  (\ref{splitn}), we have 
	$$ k_{1} + k_{2} +..........+k_{p}+.......... +k_{p^{2}} +........+k_{p^{3}}= n. $$
	Thus  the number of all  distinct  $p^3$-tuples ($k_{1}, k_{2}, k_{3},\cdots \cdots \cdots , k_{p^3}$) is $\binom{n+p^{3}-1}{p^{3}-1}$.\\

\end{flushleft}

Thus from equation (\ref{splitn}) and Theorem $\ref{theorem2.2}$ the  number of  $n$ degree representations (upto isomorphism) of a group $G$ of order $p^3$  is $\sum_{\mu=0}^{\lfloor \frac{n}{p}\rfloor} \binom{\mu +p-2}{p-2}\binom{n- \mu p+p^{2}-1}{p^{2}-1}$, when $G$ is non-abelian, whereas it is $\binom{n+p^{3}-1}{p^{3}-1}$, when $G$ is abelian.
$\hspace{5.6in}$  $\square$\\
\subsection{Degenerate invariant  spaces} \noindent From Lemma \ref{singular}, if $k_{2i} \neq k_{2i+1}$ then  all the elements of the space are degenerate Thus by the  Theorem $\ref{theorem1.1}$ and Lemma $\ref{nondeg}$, the number of $n$ degree representations whose corresponding  invariant  spaces of bilinear forms   contain only degenerate invariant bilinear forms are $\sum_{\mu=0}^{\lfloor \frac{n}{p}\rfloor} \binom{\mu +p-2}{p-2}\binom{n-\mu p+p^{2}-1}{p^{2}-1}-$ $\sum_{\ell=0}^{ \lfloor \frac{n}{2p} \rfloor}\Bigg{[}\binom{\ell+\frac{p-3}{2}}{\frac{p-3}{2}} \sum_{\lambda=0}^{\lfloor \frac{n-2p\ell}{2} \rfloor} \binom{\lambda+\frac{p^2-3}{2}}{\frac{p^2-3}{2}}\Bigg]$ in the non-abelian case and it is $\binom{n+p^{3}-1}{p^{3}-1}- \sum_{\ell=0}^{\lfloor \frac{n}{2}\rfloor}\binom{\ell+\frac{p^{3}-3}{2}}{\frac{p^{3}-3}{2}}$ in the abelian
case.\\


\noindent \textbf{Proof of  Theorem 
	$\ref{theorem1.2}$} \,

Let $A_G=\{(i,j)\,|\,\, \rho_i$ and $\rho_j$ are dual to each other\}  and  for every $(i,j)\in A_G$,  $\mathbb{W}_{(i,j)}$ = \{ $ X \in \mathbb{M}_{n}(\mathbb{F})\, | \, X_{d_{i}k_{i} \times d_jk_j}^{ij}=  C_{k_{i}\rho_{i}(g)}^{t}X_{d_{i}k_{i} \times d_jk_j}^{ij}C_{k_{j}\rho_{j}(g)}$, $\forall g \in G$ and rest blocks are zeros \}. Then for  $(i,j)\in A_G$, $\mathbb{W}_{(i,j)}$ is a  subspace of $\mathbb{M}_{n}(\mathbb{F})$. Let $X$ be an element of $\Xi_{G}'$, then 

\[C_{\rho(g)} ^{t}X C_{\rho(g)}= X \,\,and\,\, X=[ X^{ij}_{d_{i}k_{i} \times d_{j}k_{j}} ]_{(i,j) \in A_G}. \]
\noindent Existence:\\
Let $X \in \Xi_{G}' $ then for every   $ (i,j) \in A_G$,  there exists at least one  $X_{(i,j)} \in \mathbb{W}_{(i,j)}$,  such that  $\sum_{(i,j) \in A_G} X_{(i,j)} =X$.\\

\noindent Uniqueness: \\
For every  $ (i,j) \in A_G$, suppose there are  $Y_{(i,j)}$ and $X_{(i,j)} \in \mathbb{W}_{(i,j)} $, such that   $\sum_{(i,j) \in A_G}Y_{(i,j)} =X=\sum_{(i,j) \in A_G}X_{(i,j)}$, then  $\sum_{(i,j) \in A_G} X_{(i,j)} $ = $\sum_{(i,j) \in A_G} Y_{(i,j)}$ i.e.,  $Y_{(i',j')}-X_{(i',j')}$= $\sum_{\substack{(i,j) \in A_G \\(i,j) \neq (i',j')}} (X_{(i,j)}-Y_{(i,j)})$. Therefore $Y_{(i',j')}-X_{(i',j')} \in $ $\sum_{\substack{(i,j) \in A_G \\(i,j) \neq (i',j')}}  \mathbb{W}_{(i,j)}$, hence  $Y_{(i',j')}-X_{(i',j')} = O$  $\implies Y_{(i',j')}=X_{(i',j')}$ for all $(i',j') \in A_G$.\\Thus we have 
\begin{equation}
	\Xi_{G}' =\oplus_{(i,j) \in A_G} \mathbb{W}_{(i,j)} \hspace{0.1cm} and \hspace{0.1cm} dim(\Xi_{G}') = \sum_{(i,j) \in A_G} dim\mathbb{W}_{(i,j)}.
	\label{wg}
\end{equation}

\begin{flushleft}
	Now 
	as for $(i,j) \in A_G$, $\mathbb{W}_{(i,j)} = \{ X \in \mathbb{M}_{n}(\mathbb{F})\, |\,  \,  \, (i,j)^{th}$ block  '$X^{ij}$' is a  sub - matrix of order $d_{i}k_{i} \times d_jk_j $ satisfying $X^{ij}=  C_{k_{i}\rho_{i}(g)}^{t}X^{ij}C_{k_{j}\rho_{j}(g)}$, $\forall g \in G$ and rest blocks are zeros \}.
	Now by the remark  \ref{remark}, we see that for  $(i,j) \in A_G$, the sub-matrices  $X^{ij}$ in  $\mathbb{W}_{(i,j)} $  have $k_i k_j$ free variables $\&$  $\mathbb{W}_{(i,j)} \cong $ $\mathbb{M}_{k_i \times k_j }(\mathbb{F})$. Thus $\Xi_{G}'\cong $ $\oplus_{(i,j) \in A_G}\mathbb{M}_{k_i \times k_j}(\mathbb{F})$ and $dim(\mathbb{W}_{(i,j)} )= k_{i}k_j$.
\end{flushleft}


Thus  substituting these in equation  (\ref{wg}) we get the dimension of $\Xi_{G}'$. \\

\noindent \textbf{Proof of Theorem $\ref{theorem1.3}$} Follows immediately from  Lemmas  $\ref{singular}$  and  $\ref{nondeg}$ .$\hspace{2.4 in}$ $\square$\\

\section { Representation over a field of characteristic $p$.}
\begin{Remark}
	If characteristic of the field $\mathbb{F}$ is $p$ then a  group $G$ of order $p^3$ has only  trivial irreducible representation.
	Therefore  for every $g\in G$, we have
	$$\rho(g)=n\rho_{1}(g),$$
	where $\rho_{1}$ is the  trivial representation of group $G$. So  the  representation $\rho$ is a trivial  representation of degree n. i.e,  \\
	$$  \rho(g)= I_{n},\mbox{\, for  all $g \in G$.} $$
\end{Remark}
\begin{Note}
	As $X = I_n^t XI_n$,  $ \forall X \in \mathbb{M}_n(\mathbb{F})$, therefore every $X \in \mathbb{M}_n(\mathbb{F})$ gives an invariant bilinear form in the case when characteristic of the field is $p$. This is summarised in the following result.
\end{Note}
\begin{Proposition}  The space of invariant bilinear forms under an $n$ degree   representation of a group $G$ of order $p^3$ with char($\mathbb{F}$) $=p \mbox { an odd prime}$ is  isomorphic to  $\mathbb{M}_{n}(\mathbb{F})$ and  contains a  non-degenerate invariant bilinear form. 
\end{Proposition}

\vskip2mm
Thus here we have completely characterised the representations of a group of order $p^3$ each of which admits a non-degenerate invariant bilinear form over a field of characteristic different from $p$ consisting of a primitive $p^3$th root of unity. Authors hope to evaluate these results for a group of higher order in future.
\vskip2mm
\noindent {\bf Acknowledgement} The first author would like to thank UGC, India for providing the research fellowship and to the Central University of Jharkhand, India for facilitating this research work. The second author would like to express his gratitude towards Babasaheb Bhimrao Ambedkar University, Lucknow, India where he got affiliated while finalizing this paper.


\end{document}